\documentclass[leqno, fleqn]{amsart}%
\usepackage{amsmath}
\usepackage{amsfonts}
\usepackage{amssymb}%

\makeatletter \DeclareMathSymbol{\Gamma}{\mathalpha}{letters}{"00}
\DeclareMathSymbol{\Delta}{\mathalpha}{letters}{"01}
\DeclareMathSymbol{\Theta}{\mathalpha}{letters}{"02}
\DeclareMathSymbol{\Lambda}{\mathalpha}{letters}{"03}
\DeclareMathSymbol{\Xi}{\mathalpha}{letters}{"04}
\DeclareMathSymbol{\Pi}{\mathalpha}{letters}{"05}
\DeclareMathSymbol{\Sigma}{\mathalpha}{letters}{"06}
\DeclareMathSymbol{\Upsilon}{\mathalpha}{letters}{"07}
\DeclareMathSymbol{\Phi}{\mathalpha}{letters}{"08}
\DeclareMathSymbol{\Psi}{\mathalpha}{letters}{"09}
\DeclareMathSymbol{\Omega}{\mathalpha}{letters}{"0A}
\DeclareMathSymbol{\varGamma}{\mathalpha}{operators}{"00}
\DeclareMathSymbol{\varDelta}{\mathalpha}{operators}{"01}
\DeclareMathSymbol{\varTheta}{\mathalpha}{operators}{"02}
\DeclareMathSymbol{\varLambda}{\mathalpha}{operators}{"03}
\DeclareMathSymbol{\varXi}{\mathalpha}{operators}{"04}
\DeclareMathSymbol{\varPi}{\mathalpha}{operators}{"05}
\DeclareMathSymbol{\varSigma}{\mathalpha}{operators}{"06}
\DeclareMathSymbol{\varUpsilon}{\mathalpha}{operators}{"07}
\DeclareMathSymbol{\varPhi}{\mathalpha}{operators}{"08}
\DeclareMathSymbol{\varPsi}{\mathalpha}{operators}{"09}
\DeclareMathSymbol{\varOmega}{\mathalpha}{operators}{"0A}

\newcommand{\allmodesymb}[2]{\relax\ifmmode{\mathchoice
{\mbox{\fontsize{\tf@size}{\tf@size}#1{#2}}}
{\mbox{\fontsize{\tf@size}{\tf@size}#1{#2}}}
{\mbox{\fontsize{\sf@size}{\sf@size}#1{#2}}}
{\mbox{\fontsize{\ssf@size}{\ssf@size}#1{#2}}}} \else
\mbox{#1{#2}}\fi}

\makeatother
\makeatletter
\renewcommand*\subjclass[2][2000]{%
  \def\@subjclass{#2}%
  \@ifundefined{subjclassname@#1}{%
    \ClassWarning{\@classname}{Unknown edition (#1) of Mathematics%
      Subject Classification; using '2000'.}%
  }{%
    \@xp\let\@xp\subjclassname\csname subjclassname@#1\endcsname%
  }%
} \makeatother

\theoremstyle{plain}

\theoremstyle{remark}

\allowdisplaybreaks \numberwithin{equation}{section}
\begin{document}
\title[Two Problems on  Cartan Domains]
{Two Problems on  Cartan Domains}
\date{}

\author{Weiping YIN }
\thanks{Partially supported by grant from Science Research Foundation of
Academia Sinica.}

\begin{abstract}
Firstly, we consider the unitary geometry of two exceptional
Cartan domains $\Re_{V}(16)$ and $\Re_{VI}(27)$. We obtain the
explicit formulas of Bergman kernal funtion, Cauchy-Szeg\"{o}
kernel, Poinsson kernel and Bergman metric for $\Re_{V}(16)$ and
$\Re_{VI}(27)$. Secondly, we give a class of invariant
differential operators for Cartan domain $\Re$ of dimension n: If
the Bergman metric of $\Re$ is
$$ds^{2}=\sum\limits_{i,j=1}^{n}g_{ij}dz_{i}d\overline{z}_{j},
T(z,\overline{z})=(g_{ij})$$ and $$L(u)=T^{-1}(z,\overline{z})
[\frac{\partial^2u}{\partial z_i\partial\overline{z}_j}],$$then
$$L_j(u)=\{\mbox {The sum of all prinipal minors of degree}\ j\ \mbox{for}\ L(u)\}$$
is invariant under the biholomorphic mapping of $\Re$. Let $D$ be
the irreducible bounded homogeneous domain in $C^n$, $P=P(z,*)$
the Poisson kernel of $D$, then for any fixed $J(1\leq j \leq n)$
one has $L_j(P^{1/j})=0$ iff $D$ is a symmetric domain.

\end{abstract}

\maketitle

\section{Unitary Geometry on Exceptional Cartan Domains}

In 1935, E.Cartan classified all symmetric bounded domains. He
prived that there exit only six types of irreducible bounded
symmetric domains in $\mathbb{C}^n$. They can be realized as
follows:
$$\begin{array}{l}
\Re_{I}(m,n)=\{Z\in {C}^{mn}|I-Z\overline{Z}'>0,Z-(m,n)
\mbox{ matrix}\}\\
\Re_{II}{p}=\{Z\in
{C}^{p(p+1)/2}|I-Z\overline{Z}'>0,Z-\mbox{symmetric
 matric of degree}\  p\}\\
 \Re_{III}{q}=\{Z\in {C}^{q(q-1)/2}|I-Z\overline{Z}'>0,Z-\mbox{skew symmetric
 matrix of degree}\  q\}\\
 \Re_{IV}(n)=\{z=(z_1,\cdots,z_n)\in {C}^n|1+|zz'|^2-2z\overline{z}'>0,
1-|zz'|^2>0 \} \end{array}$$ or
$$=\left\{Z=\left(
\begin{array}{ccccc}
z_1&z_2&z_3&\cdots&z_{n-1}\\
 z_2& z_n & 0 &\cdots& 0 \\
 z_3&0&z_n&\cdots&0 \\
 \cdots&\cdots&\cdots&\cdots&\cdots\\
 z_{n-1}&0&0&\cdots &z_n
 \end{array}
 \right)
 \in {C}^n|\frac{1}{2\sqrt{-1}}(Z-\overline{Z}')>0\right\}.$$
 Besides the above four Cartan domains, there exist two exceptional
 Cartan domains of dimensions 16 and 27. If we denote them by $\Re_{V}(16)$ and
$\Re_{VI}(27)$ respectively, then
$$\Re_{V}(16)=\{(Z,U)\in {C}^{16}|\frac{1}{2\sqrt{-1}}(Z-\overline{Z}')
-\frac{1}{2}(U\overline{U}'+\overline{U}U')>0\},$$ where
$$Z=\left(
\begin{array}{ccccc}
z_1&z_2&z_3&\cdots&z_7\\
 z_2& z_8& 0 &\cdots& 0 \\
 z_3&0&z_8&\cdots&0 \\
 \cdots&\cdots&\cdots&\cdots&\cdots\\
 z_7&0&0&\cdots &z_8
 \end{array}
 \right),
 U=\left(
 \begin{array}{c}
 t\\uQ_1\\\vdots \\uQ_6
 \end{array}
 \right),
 \begin{array}{c}
 t=(t_1,\cdots,t_4)\in {C}^4\\
u=(u_1,\cdots,u_4)\in {C}^4
\end{array} $$
and
$$Q_1=I^{(4)},\ \
Q_2=\sqrt{-1}\left(
\begin{array}{cc}
I^{(2)}&0\\
0&-I^{(2)}
\end{array}
\right),\ \  Q_3=\left(
\begin{array}{cc}
0&I^{(2)}\\
-I^{(2)}&0
\end{array}
\right )$$
$$Q_4=\sqrt{-1}\left(
\begin{array}{cc}
0&{\begin{array}{cc}1&0\\0&-1\end{array}}\\
{\begin{array}{cc}1&0\\0&-1\end{array}}&0
\end{array}
\right),Q_5=\left(
\begin{array}{cc}
0&{\begin{array}{cc}0&1\\-1&0\end{array}}\\
{\begin{array}{cc}0&1\\-1&0\end{array}}&0
\end{array}
\right), $$$$Q_6=\sqrt{-1}\left(
\begin{array}{cc}
0&{\begin{array}{cc}0&1\\1&0\end{array}}\\
{\begin{array}{cc}0&1\\1&0\end{array}}&0
\end{array}
\right), $$
$$Q_i\overline{Q}'_j+Q_j\overline{Q}'_i=2\delta_{ij}I^{(4)},\
i,j=1,2,\cdots,6.$$
$$\Re_{VI}(27)=[(z_{11},z_{12},z_{13},z_{22},z_{23},z{33})\in
{C}^1\times{C}^8\times{C}^8\times{C}^1\times {C}^8\times{C}^1
|\frac{1}{2\sqrt{-1}}(Z-\overline{Z}')>0],$$
where$$Z=\left( \begin{array}{ccc} z_{11}&z_{12}&z_{13}\\
z'_{12}&z_{22}I^{(8)}&z_{23}\\
z'_{13}&z'_{23}&z_{33}I^{(8)} \end{array}\right),\ \
z_{23}=\left(\begin{array}{c} zT_1\\\cdots\\zT_2
\end{array}
\right),\ z=(z_1,\cdots,z_8)\in C^8 $$ and
$$T_1=\left(\begin{array}{cccc}
{\left(\begin{array}{cc}1&0\\0&-1\end{array}\right)}&&&\\
&{\left(\begin{array}{cc}-1&0\\0&-1\end{array}\right)}&&\\
&&{\left(\begin{array}{cc}-1&0\\0&-1\end{array}\right)}&\\
&&&{\left(\begin{array}{cc}-1&0\\0&-1\end{array}\right)}
\end{array}
\right),$$
$$T_2=\left(\begin{array}{cccc}
{\left(\begin{array}{cc}0&1\\1&0\end{array}\right)}&&&\\
&{\left(\begin{array}{cc}0&-1\\1&0\end{array}\right)}&&\\
&&{\left(\begin{array}{cc}0&-1\\1&0\end{array}\right)}&\\
&&&{\left(\begin{array}{cc}0&1\\-1&0\end{array}\right)}
\end{array}
\right),$$
$$T_3=\left(\begin{array}{cccc}
&{\left(\begin{array}{cc}1&0\\0&1\end{array}\right)}&&\\
{\left(\begin{array}{cc}1&0\\0&-1\end{array}\right)}&&&\\
&&&{\left(\begin{array}{cc}-1&0\\0&-1\end{array}\right)}\\
&&{\left(\begin{array}{cc}1&0\\0&1\end{array}\right)}&
\end{array}
\right),$$
$$T_4=\left(\begin{array}{cccc}
&{\left(\begin{array}{cc}0&1\\-1&0\end{array}\right)}&&\\
{\left(\begin{array}{cc}0&1\\1&0\end{array}\right)}&&&\\
&&&{\left(\begin{array}{cc}0&-1\\1&0\end{array}\right)}\\
&&{\left(\begin{array}{cc}0&-1\\1&0\end{array}\right)}&
\end{array}
\right),$$
$$T_5=\left(\begin{array}{cccc}
&&{\left(\begin{array}{cc}1&0\\0&1\end{array}\right)}&\\
&&&{\left(\begin{array}{cc}1&0\\0&1\end{array}\right)}\\
{\left(\begin{array}{cc}1&0\\0&-1\end{array}\right)}&&&\\
&{\left(\begin{array}{cc}-1&0\\0&-1\end{array}\right)}&&
\end{array}
\right),$$

$$T_6=\left(\begin{array}{cccc}
&&{\left(\begin{array}{cc}0&1\\-1&0\end{array}\right)}&\\
&&&{\left(\begin{array}{cc}0&1\\-1&0\end{array}\right)}\\
{\left(\begin{array}{cc}0&1\\1&0\end{array}\right)}&&&\\
&{\left(\begin{array}{cc}0&1\\-1&0\end{array}\right)}&&
\end{array}
\right),$$
$$T_7=\left(\begin{array}{cccc}
&&&{\left(\begin{array}{cc}1&0\\0&-1\end{array}\right)}\\
&&{\left(\begin{array}{cc}-1&0\\0&1\end{array}\right)}&\\
&{\left(\begin{array}{cc}1&0\\0&-1\end{array}\right)}&&\\
{\left(\begin{array}{cc}1&0\\0&1\end{array}\right)}&&&
\end{array}
\right),$$
$$T_8=\left(\begin{array}{cccc}
&&&{\left(\begin{array}{cc}0&1\\1&0\end{array}\right)}\\
&&{\left(\begin{array}{cc}0&-1\\-1&0\end{array}\right)}&\\
&{\left(\begin{array}{cc}0&1\\1&0\end{array}\right)}&&\\
{\left(\begin{array}{cc}0&-1\\1&0\end{array}\right)}&&&
\end{array}
\right).$$
$$T_iT'_j+T_jT'_i=2\delta_{ij}I^{(8)},i,j=1,2,\cdots,8.$$
For the first four types of Cartan domains, Hua and Lu obtained
many results[2,3].Now we consider the unitary geometry on the
exceptional Cartan domains $\Re_{V}(16)$ and $\Re_{VI}(27)$.\\
 {\bf I Bergman kernel function}\\
 The mapping
 $$\begin{array}{l}
 W=A[Z-\frac{1}{2}(Z_0+\overline{Z}'_0)-\sqrt{-1}(U\overline{U}'_0+
 \overline{U}_0U')+\frac{\sqrt{-1}}{2}(U_0\overline{U}'_0+\overline{U}_0
 U'_0)]A'\\
 R=A(U-U_0)
 \end{array}$$
 is a holomorphic automorphism of  $\Re_{V}(16)$ ,where
 $$R=\left(\begin{array}{c}r\\sQ_1\\\vdots\\sQ_6\end{array}\right),\begin{array}{c}
 r=(r_1,\cdots,r_4)\in {C}^{4}\\ s=(s_1,\cdots,s_4)\in {C}^{4}
 \end{array}$$
 $$A=\left(\begin{array}{cccc}
a_{11}& a_{12}&\cdots&a_{17}\\
&a_{22}&&\\
&&\ddots&\\
&&&a_{22}
\end{array}\right)\in {R}^8$$
which maps point $(Z_0,U_0)(\in \Re_{V}(16))$ into point
$(\sqrt{-1}I,0)$, and
$$(A'A)^{-1}=\frac{Z_0-\overline{Z}'_0}{2\sqrt{-1}}-\frac{1}{2}(U_0\overline{U}'_0
+\overline{U}_0U'_0),\
a_{22}^{-2}=\frac{z_8^0-\overline{z}_8^0}{2\sqrt{-1}}-u_0\overline{u}'_0.$$
By direct calculations, we have
$$\frac{\partial{(W,R)}}{\partial{(Z,U)}}=
\left(\begin{array}{ccccc} a_{11}^2&&&&\\
&a_{11}a_{22}I^{(6)}&&&\\
&&a_{22}^2&&\\
&&&a_{11}I^{(4)}&\\
&&&&a_{22}I^{(4)} \end{array}\right)$$ and
$$\det\left[\frac{\partial{(W,R)}}{\partial{(Z,U)}}
\overline{\frac{\partial{(W,R)}}{\partial{(Z,U)}}'}\right] =
(a_{11}a_{12})^{24}=(a_{11}^2a_{22}^2)^{12}/a_{22}^{120} = \det
(A'A)^{12}/a_{22}^{120}.$$ Hence, the Bergman kernel function of
$\Re_{V}(16)$ is given by(up to a constant factor)
\begin{eqnarray*}
K_{V}(Z,U,Z,U)&=&\frac{\left[\frac{1}{2\sqrt{-1}}(z_8-\overline{z}_8)-u\overline{u}'
\right]^{60}}{\det\left[(2\sqrt{-1})^{-1}(Z-\overline{Z}')-\frac{1}{2}(U\overline{U}'
+\overline{U}U')\right]^{12}}\\
&=&\{(2\sqrt{-1})^{-1}(z_1-\overline{z}_1)((2\sqrt{-1})^{-1}(z_8-\overline{z}_8))\\&&
-\sum_{j=1}^{6}[(2\sqrt{-1})^{-1}(z_{j+1}-\overline{z}_{j+1})-
\frac{1}{2}(uQ_j\overline{t}'+t\overline{Q}'_j\overline{u}')]^2\}^{-12}.
\end{eqnarray*}

Mapping $$W=A\left[Z-\frac{1}{2}(Z_0+\overline{Z}'_0)\right]A'$$
is a holomorphic automorphism of $\Re_{VI}(27)$, where
$$A=\left(\begin{array}{ccc}
a_{11}&a_{12}&a_{13}\\
0&a_{22}I^{(8)}&a_{23}\\
0&0&a_{33}I^{(8)} \end{array} \right),
a_{23}=\left(\begin{array}{c}aT_1\\aT_2\\\vdots\\aT_8
\end{array}\right),a=(a_1,\cdots,a_8)\in {R}^{8}.$$
which maps point$Z_0$($\in R_{VI}$(27)) into point $\sqrt{-1}$I
and$$
(A'A)^{-1}=(2\sqrt{-1}){-1}(Z_0-\overline{Z_0}'),(a_{22}a_{33})^{-2}=(\frac{{z_{22}}^0-
{\overline{z_{22}}}^0}{2\sqrt{-1}})(\frac{{z_{33}}^0-
{\overline{z_{33}}}^0}{2\sqrt{-1}})-(\frac{z_0-
\overline{z_0}}{2\sqrt{-1}})\overline{(\frac{z_0-
\overline{z_0}}{2\sqrt{-1}})}'.
$$

By direct calculations, we have

$$ \left(
\begin{array}{cccccc}
a^2_{11}&&&&&\\&a_{11}a_{22}I^{(8)}&&&&\\&&a_{11}a_{33}I^{(8)}&&&\\&&&{a^2_{22}}\\&\ast&&&{a_{22}a_{33}I^{(8)}}&\\
&&&&&{a^2_{33}}
\end{array}
\right)=\frac{\partial{W}}{\partial{Z}}.
 $$
and
$$
\det(\frac{\partial{W}\overline{\partial{W}}}{\partial{Z}\partial{Z}})
=(a_{11}a_{22}a_{33})^{36}=(a^2_{11}a^{16}_{22}a^{16}_{33})^{18}/(a^2_{22}a^2_{33})^{126}
$$
$$
=(\det A'A)^{18}/(a^2_{22}a^2_{33})^{126}
$$
So, we obtain the Bergman kernel function of $R_{VI}(27)$ as
follows(up to a constant factor):
$$
K_{VI}(Z,\overline{Z})=\frac{[\frac{1}{2\sqrt{-1}}(z_{22}-
\overline{z_{22}})\frac{1}{2\sqrt{-1}}(z_{33}-
\overline{z_{33}})-\frac{1}{2\sqrt{-1}}(z-
\overline{z})\overline{\frac{1}{2\sqrt{-1}}(z-
\overline{z})}']^{126}}{\det [\frac{1}{2\sqrt{-1}}(z-
\overline{z})]^{18}}.
$$
{\bf 2.Cauchy-Szeg\"{} Kernels and Poisson Kernels}

If the point $(X,V)$ belongs to Solov boundary of $R_{V}(16)$,
then the $(X,V)$ satisfies the following equation:
$$
\frac{X-\overline{X}'}{2\sqrt{-1}}=\frac{1}{2}(V\overline{V}'+\overline{V}V'),
V=\left(
\begin{array}{cccc}
u\\ vQ_{1}\\\vdots\\ vQ_{6}
\end{array}
\right), u=(u_1,\cdots,u_4)\in C^4, v=(v_1,\cdots,v_4)\in C^4.
$$
By direct calculation, the Cauchy-Szeg\"{o} kernel of $R_{V}(16)$
is
$$
H_V(Z,U;X,V)=\frac{[(z_8-\overline{x_8})/(2\sqrt{-1})-u\overline{v}']^{30}}
{\det[(Z-\overline{X}')/(2\sqrt{-1})-\frac{1}{2}(U\overline{V}'+\overline{U}V')]^6}
$$
(up to a constant factor).

And the Poisson kernel of $R_{V}(16)$ is given by (up to a
constant factor)
$$
P_V(Z,U;X,V)
=\frac{\det[(Z-\overline{Z}')/(2\sqrt{-1})-(U\overline{U}'+\overline{U}U')/2]^6|
(z_8-\overline{x_8})/(2\sqrt{-1})-u\overline{v}'|^{60}}
{[(z_8-\overline{z_8})/(2\sqrt{-1})-u\overline{u}']^{30}
|\det[(Z-\overline{X}')/(2\sqrt{-1})-(U\overline{V}'+\overline{U}V')/2]|^{12}}
$$
If X belong to the silov kernel boundary of $R_{VI}(27)$, then we
have
$$X=\left[
\begin{array}{ccc}
X'_{11}&X_{12}&X_{13}\\
X'_{12}&X_{22}I^{(8)}&X_{23}\\X_{13}&X'_{23}&X_{33}I^{(8)}
\end{array}
\right], X_{23}=\left[
\begin{array}{c}
XT_1\\ \vdots\\XT_8
\end{array}
\right], X=(X_1\cdots,X_8)\in \mathbf{R^8}
 $$
By direct calculations, the Cauchy-Szeg\"{o} kernel of
$R_{VI}(27)$ is given by (up to a constant factor)
$$
H_{VI}(Z,X)=\frac{(z_{22}- \overline{x_{22}})/(2\sqrt{-1})(z_{33}-
\overline{z_{33}})/(2\sqrt{-1})-(z-
\overline{x})/(2\sqrt{-1})[\overline{(z-
\overline{x})/(2\sqrt{-1})']}^{63}}{[{\det(Z-
\overline{X}')/(2\sqrt{-1})}]^{9}}.
$$
And the Poisson kernel of $R_{VI}(27)$ is (up to a constant
factor)
$$
P_{VI}(Z,X)=\frac{\det
(\frac{Z-\overline{Z}'}{2\sqrt{-1}})^9|(\frac{{z_{22}}-
{\overline{z_{22}}}}{2\sqrt{-1}})(\frac{{z_{33}}-
{\overline{z_{33}}}}{2\sqrt{-1}})-(\frac{z-
\overline{x}}{2\sqrt{-1}})\overline{(\frac{z-
\overline{x}}{2\sqrt{-1}})}'|^{126}}{|\det
(\frac{Z-\overline{X}'}{2\sqrt{-1}})|^{18}[(\frac{z_{22}-
\overline{z_{22}}}{2\sqrt{-1}})(\frac{{z_{33}}-
{\overline{z_{33}}}}{2\sqrt{-1}})-(\frac{z-
\overline{z}}{2\sqrt{-1}})\overline{(\frac{z-
\overline{z}}{2\sqrt{-1}})}']^{63}}
$$

\section{A Class of Invariant Differential Operators on Cartan Domains and Their solutions}

1.We consider the Cartan domain of first type
\[
R_{I}(m,n)=(Z|I-Z\overline{Z'}>0,Z-(m,n)\ matrix).
\]
It is well known that the Bergman kernel function is given by (up
to a consant factor)
\[
K_{I}(Z,\overline{Z})=[det(I-Z\overline{Z'})]^{-(m+n)}.
\]
Let
\[
(g_{j\alpha},\overline{k\beta})=(\frac{\partial^{2}lgK_{I}(Z,\overline{Z})}{\partial
z_{j\alpha}\partial {\overline{z}}_{k\beta}}).
\]
Then $(g_{j\alpha},\overline{k\beta})=(I-Z\overline{Z'})^{-1}\cdot
X(I-Z\overline{Z'})^{-1}=T_{I}(Z,\overline{Z})=$Bergman metric
matrix of $R_{I}(m,n)$,
$ds_{I}^{2}=\sum^{m}_{j,k=1}\sum^{n}_{\alpha,\beta=1}g_{j\alpha,\overline{k\beta}}dz_{j\alpha}d\overline{z}_{k\beta}=$
Bergman metric matrix of $R_{I}$(up to a const.factor).And
\[
det(g_{j\alpha},\overline{k\beta})=[det(I-Z\overline{Z'})]^{-(m+n)}=K_{I}(Z,\overline{Z}).
\]
Let
$$\displaystyle L=T^{-1}_{1}(Z,\overline{Z})\frac{\partial^{2}}{\partial z'\partial \overline{z}},\ \ \ \ \ \
L(u)=T^{-1}_{1}(Z,\overline{Z})\frac{\partial^{2}u}{\partial
z'\partial \overline{z}},$$

where
$$\frac{\partial}{\partial z}
=\left(\frac{\partial}{\partial z_{11}},\frac{\partial}{\partial
z_{12}},...,\frac{\partial}{\partial
z_{1n}},\frac{\partial}{\partial
z_{21}},...,\frac{\partial}{\partial
z_{2n}},...,\frac{\partial}{\partial
z_{m1}},...,\frac{\partial}{\partial z_{mn}}\right),$$

$$\frac{\partial^{2}}{\partial z'\partial {\overline{z}}}
=\left(\frac{\partial}{\partial
z}\right)'\left(\overline{\frac{\partial}{\partial z}}\right)
=\left(
\begin{array}{ccc}
\frac{\partial^{2}}{\partial z_{11}\partial
\overline{z_{11}}}&\cdots\cdots\cdots&
\frac{\partial^{2}}{\partial z_{11}\partial \overline{z_{mn}}}\\\cdots&\cdots\cdots\cdots&\cdots\\
\frac{\partial^{2}}{\partial z_{mn}\partial
\overline{z_{11}}}&\cdots\cdots\cdots&
\frac{\partial^{2}}{\partial z_{mn}\partial \overline{z_{mn}}}
\end{array}
\right).$$

If $W=f(Z)$ is the holomorphic automorphism of $R_{1}(m,n)$, then

$$T^{-1}_{1}(Z,\overline{Z})\frac{\partial^{2}u(Z)}{\partial z'\partial{\overline{z}}}
=\left(\overline{\frac{\partial W}{\partial
Z}}\right)'^{-1}T^{-1}_{1}(W,\overline{W})T^{-1}_{1}(Z,\overline{Z})\frac{\partial^{2}u(f^{-1}(W))}{\partial
W'\partial \overline{W}}{\left(\overline{\frac{\partial
W}{\partial Z}}\right)}'\ \ \ \ \ \ (4)$$

Let
$$L_{j}(u)=\rm \{ The\  sum \ of \ all \ principal \ minors \ of \ degree \ j \ for \ L(u)\}$$
then $L_{j}(u)$ is an invariant differential operator of
$R_{1}(m,n)(j=1,2,...,mn)$.

In fact, from $(4)$ we know that
$$T^{-1}_{1}(Z,\overline{Z})\frac{\partial^{2} u}{\partial z'\partial \overline{z}} \rm \ is \ similar \ to \ T^{-1}_{1}(W,\overline{W})\frac{\partial^{2} u(f^{-1}(W))}{\partial W'\partial \overline{W}}$$
Suppose
$$F(\lambda)
=det\left[\lambda I-T^{-1}_{1}(Z,\overline{Z})\frac{\partial^{2}
u}{\partial z'\partial \overline{z}}\right]$$ is the
characteristic polynomial for $\displaystyle
T^{-1}_{1}\frac{\partial^{2} u(Z)}{\partial z'\partial
\overline{z}},$ then the coefficient of $\lambda^{mn-j}$ is the
$L_{j}(u)(j=1,2,......,mn)$\ \ (\ up \ to \ sign$\pm$). But the
similar matriceshave the same characteristic polynomial .So
$L_{j}(u)$ is an invariant differential operator.

Specifically,
$$L_{1}(u)=tr\left[T^{-1}_{1}(Z,\overline{Z})\frac{\partial^{2} u}{\partial z'\partial \overline{z}}\right]$$
is the Laplace-Beltrami operator of $(R_{1},ds^{2}_{1})$.

And
$$L_{mn}(u)=\rm tr\left[T^{-1}_{1}(Z,\overline{Z})\frac{\partial^{2} u}{\partial z'\partial \overline{z}}\right]$$
is the complex Monge-Ampere operator.

2. Let
$$\displaystyle P_{j}(Z,U)=\frac{\rm det(I-Z\overline{Z}')^{n/j}}{|\rm det(I-Z\overline{U})'|^{2n/j}}\ \ \ (j=1,2,......,mn),$$
$$U\overline{U}'=I,U \ \rm is \ a \ (m,n) \ matrix.$$
Then
$$L_{j}(p_{j}(Z,U))=0,\ \ \ j=1,2,......,mn.$$
In fact,if we expand the $P_{j}(Z,U)$ around the point $Z=0$:
$$\displaystyle P_{j}(Z,U)
=\left[1-\frac{n}{j}\rm tr(Z\overline{Z}')+...\right]
\left[1+\frac{n}{j}\rm
tr(Z\overline{U}')+...\right]\left[1+\frac{n}{j}\rm
tr(U\overline{Z}')+...\right],$$ Then,we have
$$\displaystyle\frac{\partial^{2}P_{j}(Z,U)}{\partial z_{j\alpha}\partial{\overline{z}_{k\beta}}}|_{Z=0}
=\displaystyle\left\{
\begin{array}{ll}
\displaystyle -\frac{n}{j}+\left(\frac{n}{j}\right)^{2}|u_{j\alpha}|^{2},&\quad \rm if (j,\alpha)=(k,\beta)\\
\displaystyle\left(\frac{n}{j}\right)^{2}\overline{u}_{j\alpha}u_{k\beta},&\quad
\rm if (j,\alpha)\neq (k,\beta),
\end{array}
\right.$$

where
$$U=(u_{j\alpha})\ \ \ \ \ \rm and \ U\overline{U}'=I^{(m)}.$$
So
$$L(P_{j})_{Z=0}=-\displaystyle\left(\frac{n}{j}\right)I^{(mn)}+\left(\frac{n}{j}\right)^{2}\overline{u}'u,$$
where
$$u=(u_{11}...u_{1n}u_{21}...u_{2n}...u_{m1}...u_{mn}),\ \ \ u\overline{u}'=m.$$
Then we have
$$\displaystyle L_{j}(P_{j})|_{Z=0}
=\left[C^{j}_{mn}\left(1-\frac{n}{j}u\overline{u}'\right)
+C^{j}_{mn}\frac{(n)(mn-j)}{jmn}m\right]\left(-\frac{n}{j}\right)^{j}=0,$$
where $C^{j}_{mn}=\frac{(mn)!}{j!(mn-j)!}.$

Mapping
$$W=(AZ+B)(CZ+D)^{-1}=(\overline{A}'+Z\overline{B}')^{-1}(\overline{C}'+Z\overline{D}')$$
is a holomorphic automorphism of $R_{1}(m,n)$ ,which maps $U$ into
$V$ and maps point $Z_{0}=-A^{-1}B$ into point $W=0$
 and we have
 $$I-Z_{0}\overline{Z}_{0}'=(\overline{A}'A)^{-1},\ \ \ I-\overline{Z_{0}}'Z_{0}=(\overline{D}'D)^{-1},
 \ \ \ \overline{C}'\overline{D}'^{-1}=A^{-1}B.$$
 Then, by direct calculations, we have
 $$P_{j}(W,V)=\frac{P_{j}(Z,U)}{P_{j}(Z_{0},U)}.\ \ \ \ \ \ \ (j=1,2,......,mn)$$
So we have
$$L_j(P_j(Z,U))|_{Z=Z_0}=P_j(Z_0,U)L_j(P_j(W,V))|_{W=0}=0,
(j=1,2,\dots, mn)$$
 This is
 $$L_{j}(P_{j}(Z,U))=0,\ \ \ \ \ \ \ j=1,2,......,mn.$$
 If $f(U)$ is continuous on $SR_{1}(m,n)$ (the Silov boundary of $R_{1}(m,n)$),

 $$SR_{1}(m,n)=\{U^{m,n}|U\overline{U}'=I^{(m)}\},$$and\ \ let$$Y_{j}(Z)=\displaystyle \int_{SR_{1}(m,n)}f(U)P_{j}(Z,U)\dot{U},\ \ \ \ \ j=1,2,...,mn,$$

then \ \ \ we\ \ have$$L_{1}(Y_{1}(Z))=0,$$ For the other Cartan
domain, such as $R_{II}, R_{III}, R_{IV}, R_{V}$ and $ R_{VI}$, we
have similar results.

3. According  to the above idea, we can obtain some other
invariant differential operators and some holomorphic invariants
for bounded domains in $C^{n}$. If $D_{1}$ be the bounded domain
in $C^{n}$, and $w=f(z)$ be the biholomorphic mapping of $D_{1}$,
which maps $D_{1}$ onto domain $D_{2}$ in $C^{n}$. let $T_{i}$ be
the Bergman  metric of $D_{i}(i=1,2)$, and
   $$R_{1}=-\displaystyle (\frac{\partial^{2}log det T_{1} }{\partial z_{i}\partial \overline z_{j}}),$$
   $$R_{2}=-\displaystyle (\frac{\partial^{2}log det T_{2} }{\partial w_{i}\partial \overline w_{j}}).$$
Then we have
   $$1^{\circ}.  \ \ T ^{-1}_{1}R_{1}=\displaystyle  \left(\overline{\frac{\partial W}{\partial Z}}\right)'^{-1}T^{-1}_{2}R_{2} \left(\overline{\frac{\partial W}{\partial Z}}\right)', $$
   $$R_{1}T ^{-1}_{1}=\displaystyle  \left(\frac{\partial W}{\partial Z}\right)R_{2}T^{-1}_{2}\left(\frac{\partial W}{\partial Z}\right)^{-1},$$

  $$ 2^{\circ}.\ \  T ^{-1}_{1}\overline d(dT_{1}\cdot T ^{-1}_{1})T_{1} =\displaystyle \left(\overline{\frac{\partial W}{\partial Z}}\right)^{-1}T ^{-1}_{2}\overline d(dT_{2}\cdot T ^{-1}_{2})T_{2}
  \left(\overline{\frac{\partial W}{\partial Z}}\right)'.$$
Hence
   $$\triangle_{j}(R)=\rm \{The \ sum \ of \ all\  j\  by \ j \ principal\  minors \ of \ T^{-1}_{1}R_{1}\},$$
   $$\overline \triangle_{j}(R)=\rm \{The \ sum \ of \ all\  j\  by \ j \ principal\  minors \ of \ R_{1}T^{-1}_{1}\},$$
   $$\triangle_{j}(T)=\rm \{The \ sum \ of \ all\  j\  by \ j \ principal\  minors \ of \overline d(dT_{1}\cdot T ^{-1}_{1})   \} $$
are \ \  invariant\ \  under\ \  the\ \  biholomorphic\ \ mapping\
\ of\ \ $D_{1}$. $(j=1,2,......,n)$  Also\ \,we\ \ have\ \
    that

    $$\triangle^{N}_{j}(R)=\rm \{The \ sum \ of \ all\  j\  by \ j \ principal\  minors \ of \ (T^{-1}_{1}R_{1})^{N}\}, $$
    $$\overline \triangle^{N}_{j}(R)=\rm \{The \ sum \ of \ all\  j\  by \ j \ principal\  minors \ of \ (R_{1}T^{-1}_{1})^{N}\}, $$
    $$\triangle^{N}_{j}(T)=\rm \{The \ sum \ of \ all\  j\  by \ j \ principal\  minors \ of \ (\overline d(dT_{1}\cdot T ^{-1}_{1}))^{N} \}, $$
    $$\triangle^{r}_{j}(L)=\rm \displaystyle \{The \ sum \ of \ all\  j\  by \ j \ principal\  minors \ of \ (T^{-1}_{1} \frac{\partial^{2}u}{\partial z^{\prime}\partial \overline z})^{N}\}, $$
    $$\overline \triangle^{N}_{l}(L)=\rm \displaystyle \{The \ sum \ of \ all\  j\  by \ j \ principal\  minors \ of \ ( \frac{\partial^{2}u}{\partial z^{\prime}\partial \overline z})^{N}T^{-1}_{1}\}, $$
are  invariant under the biholomorphic  mapping  of  $D_{1}$,

where   $N$ are positive  integers .

If   $D_{1}$\ \rm be\ \   the\ \   homogeneous\ \   domain ,\ \
then

 $ -T^{-1}_{1}R_{1}=I^{(n)}$ ,

 so $ \triangle_{j}(R)=\overline \triangle_{j}(R)=+C^{j}_{n}(-1)^{j}$.If  $D_{1}$\ \  be\ \  the\ \  irreducible\ \  Cartan\ \  domain\ \  , then\

   $\overline d(dT_{1}\cdot T ^{-1}_{1}|_{z=0}\overline ddT_{1}\cdot T ^{-1}_{1}$,
so  we  have $$\triangle_{1}(T)=\rm  Bergman \ \  metric\ \  of\ \
D_{1}$$.

4. \ Suppose  D  is  a  bounded  homogeneous  domain  in ${C}^{n}$
\  containing \ \  the \ \  origin.

Let  $T(Z,\overline Z),  K(Z,\overline Z),  P(Z,U)\rm$  be\ \ the\
\   Bergman \ \  metric\ \ matrix.  Bergman

kernel   function  ,formal  Poission  kernel  of  D  respectively.

By  a  linear  mapping  and  from  the  relations  between
$T(Z,\overline Z),  K(Z,\overline Z)$  and  $P(Z,U)$,

 we  can  assume  that

$$T|_{Z=0}=\lambda I,     \displaystyle \left(\frac{\partial^{2}logP}{\partial z^{\prime}\partial\overline z}\right)|_{z=0}=-\mu I,  \lambda \neq 0,   \mu \neq 0.$$

Let

$$\rm L_{j}(\mu)=\left \{the\  sum\  of\  all\  princial\   minor\   of\   degree\   j\   for\   L(\mu)\ \right\}$$

$$\displaystyle P_{j}(Z,U)=P(Z,U)^{\frac{1}{j}}$$

Then

$$L_{j}(P_{j}(Z,U))=0$$ if and only if $$L_{1}(P(Z,U))=0,\ \ \ \ \ \ \ \ (\ast)\ \ \ \ $$

where

$$j=1,2,...,n.$$

{\bf{Proof }} We  have  $$L(P_{j}(Z,U))= P_{j}(Z,U)G(
P_{j}(Z,U))$$ where $G(u)=T^{-1}(Z,\overline Z)\displaystyle
\left[(\frac{\partial ^{2}logu}{\partial
z^{\prime}\partial\overline z})+ (\frac{\partial logu}{\partial
z^{\prime}})(\frac{\partial logu}{\partial \overline z})\right]$

Let$$G_{j}(u)=\rm \{the\  sum\  of\  all\  principal\  minors\
of\  degree\  j\  for\  G(u)\},$$

then $ G_{j}(u)$  is\   also\  an\   invariant\   differential\
operators\ .(j=1,2,...,n).

If\  $W=f(Z)$\  be\  the biholomophic\  automrophism\  of\  $D$,\
which\  maps\  $Z_{0}$\  into\  O \ and

 $f(U)=V$, then  $P(Z,U)=P(Z_{0},U)P(W,V)$.so it is sufficient to

  prove that $(\ast)$ holds for Z=0.

  Because
  $$G(P_{i}(Z,U))=T^{-1}\left[\frac{1}{j}(\frac{\partial^{2}logP}{\partial z^{\prime}\partial \overline z})+\frac{1}{j^{2}}(\frac{\partial logP}{\partial z^{\prime}})(\frac{\partial logP}{\partial \overline z})\right],$$

  let$$\frac{\partial logP}{\partial z}|_{z=0}=(\alpha_{1},\alpha_{2},...,\alpha_{n})=\alpha.$$

Then we have
$$G(P_{j})|_{z=0}=-\left(\frac{\lambda\mu}{j}\right)I+\displaystyle
\frac{1}{j^{2}}\alpha^{\prime}\overline \alpha\ \ \
(let(\lambda\mu)^{-1}=b,\alpha \overline \alpha =a)$$

$$ =\left(\frac{-\lambda\mu}{j}\right)[I-(\frac{b}{j})\alpha^{\prime}\overline \alpha]$$
There exists an unitary matrix H such that
$$
I-(b/j)\alpha'\overline{\alpha}=H(I-(b/j)\left(\begin{array}{cc}
\alpha&0\\
 0&0
\end{array}
\right))H^{-1},
$$
So,
$$
L_j(P_j)|_{Z=0}=0\;iff \;X_j=0,
$$
Where
$$
X_j=\{the\; sum \;of \;all \;principal \;minors \;of \;degree\; j
\;for [I-\left(\begin{array}{cc}
\alpha(b/j)&0\\
 0&0
\end{array}
\right)]\}.
$$
But
$$
X_j=C^{j-i}_{n-1}(1-(\alpha
b/j))+{C^{j}_{n-1}}={C^{j}_{n}}(n-ab)/n.
$$
So
$$
X_1=0 \;iff \;n-ab=0.
$$
But
$$
X_1=n-ab
$$
and
$$
X_1=0\;iff\;C_1(P_1)|_{Z=0}=0\;iff \; L_1(P)|_{Z=0}=0
$$
So we complete the proof.
\\
From this we have the following theorem.
\\
{\bf{Theorem}}: If $D$ be the irreducible bounded homogenous
domain in $C^n$, then for any fixed $j(1\leq j\leq n)$ one has
$$
L_j(P_j)=0$$ if and only if $D$ is a symmetric domain.

{\bf{Proof}}: If $D$ be the irreducible symmetric domain, then
from II.2, we have $L_j(P_j)=0$.

If $L_j(P_j)=0$, then $L_1(P)=0$, i.e. the Possion kernel of D is
annihilated by the Laplace- Beltrami operator of $D$ under the
Bergman metric. Then from the Theorem 8 of [1], that $D$ must be
symmetric.

\end{document}